\documentclass[12pt]{amsart}

\usepackage{times,epsf,amssymb,amsmath,rlepsf,hyperref}

\begin{document}

\newtheorem{thm}{Theorem}[section]
\newtheorem{lem}[thm]{Lemma}
\newtheorem{cor}[thm]{Corollary}

\theoremstyle{definition}
\newtheorem{defn}{Definition}[section]

\theoremstyle{remark}
\newtheorem{rmk}{Remark}[section]

\def\square{\hfill${\vcenter{\vbox{\hrule height.4pt \hbox{\vrule
width.4pt height7pt \kern7pt \vrule width.4pt} \hrule height.4pt}}}$}

\def\T{\mathcal T}

\newenvironment{pf}{{\it Proof:}\quad}{\square \vskip 12pt}

\title{Embedded Plateau Problem}
\author{Baris Coskunuzer}
\address{Koc University \\ Department of Mathematics \\ Sariyer, Istanbul 34450 Turkey}
\email{bcoskunuzer@ku.edu.tr}
\thanks{The author is partially supported by EU-FP7 Grant IRG-226062 and TUBITAK Grant 109T685}
%\subjclass{Primary 53A10, Secondary 57M50}

\maketitle

% User definitions:

\newcommand{\SH}{S^2_{\infty}(\mathbf{H}^3)}
\newcommand{\PI}{\partial_{\infty}}
\newcommand{\SI}{S^2_{\infty}}
\newcommand{\BHH}{\mathbf{H}^3}
\newcommand{\CH}{\mathcal{C}(\Gamma)}
\newcommand{\BH}{\mathbf{H}}
\newcommand{\BR}{\mathbf{R}}
\newcommand{\BC}{\mathbf{C}}
\newcommand{\BZ}{\mathbf{Z}}
\newcommand{\BN}{\mathbf{N}}

\begin{abstract}

We show that if $\Gamma$ is a simple closed curve bounding an embedded disk in a closed $3$-manifold $M$, then there
exists a disk $\Sigma$ in $M$ with boundary $\Gamma$ such that $\Sigma$ minimizes the area among the embedded disks
with boundary $\Gamma$. Moreover, $\Sigma$ is smooth, minimal and embedded everywhere except where the boundary
$\Gamma$ meets the interior of $\Sigma$. The same result is also valid for homogenously regular manifolds with
sufficiently convex boundary.

\end{abstract}

\section{Introduction}

The Plateau problem asks the existence of an area minimizing disk for a given simple closed curve in a manifold $M$.
This problem was solved for $\BR^3$ by Douglas \cite{Do}, and Rado \cite{Ra} in the early 1930s. Later, it was
generalized by Morrey for Riemannian manifolds \cite{Mo}. Then, regularity (nonexistence of branch points) of these
solutions was shown by Osserman \cite{Os}, Gulliver \cite{Gu} and Alt \cite{Al}. However, these area minimizing disks
may not be embedded, even though the curves bound an embedded disk in the ambient manifold. They might have self
intersections (See Figure 1).

In the following decades, the question of embeddedness of the area minimizing disk was studied: For which curves are
the area minimizing disks  embedded? The first such condition ensuring the embeddedness of the disk was due to Rado. In
the early 1930s, he showed that if the curve can be projected onto a convex curve in a plane, then it bounds a unique
embedded minimal disk which is a graph over the plane. Osserman conjectured that if the curve is extreme (lies in the
boundary of its convex hull), then the area minimizing disk spanning the curve must be embedded. In the late 1970s,
Gulliver and Spruck proved that if the total curvature of an extreme curve is less than $4\pi$ then the solution to the
Plateau problem is embedded \cite{GS}. Later, Almgren-Simon \cite{AS} and Tomi-Tromba \cite{TT} showed the existence of
an embedded solution for extreme curves. Then, Meeks and Yau proved the Osserman's conjecture in full generality: Any
solution to the Plateau problem for an extreme curve must be embedded \cite{MY1}. Recently, Ekholm, White, and
Wienholtz generalized Gulliver-Spruck's embeddedness result by removing extremeness condition from the curves
\cite{EWW}. Also, recently, Hass, Lagarias and Thurston \cite{HLT} gave interesting results about the isoperimetric
inequalities for embedded disks in $\BR^3$.

On the other hand, a different version of the Plateau problem was studied after 1960s. This version asks the existence
of area minimizing surface for a given simple closed curve. If there is no restriction on the topological type of the
surface, Geometric Measure Theory gives a positive solution for this question. Federer et al solved the problem and
showed the existence of a surface which minimizes area among all surfaces with the given boundary \cite{Fe}. Moreover,
any such surface must be embedded for any simple closed curve.

If we come back to the disk case, there is a relevant result about the same question due to Meeks-Yau \cite{MY3}. They
give a necessary condition for a sufficiently smooth simple closed curve in a $3$-manifold to bound a embedded minimal
disk. In particular, they showed that for a sufficiently smooth simple closed curve $\Gamma$ in a $3$-manifold $M$, in
order to bound a strictly stable embedded minimal disk in $M$, $\Gamma$ must be an extreme curve in some sense (See
Theorem 4.1).

In this paper, we are approaching to the embeddedness question from a different direction. Instead of considering the
question that {\em ``for which curves must the area minimizing disks be embedded?"}, we analyze the structure of the
surface which minimizes area among the embedded disks whose boundary is any given simple closed curve.

\vspace{0.3cm}

\noindent {\bf Embedded Plateau Problem:} Let $\Gamma$ be a simple closed curve in a manifold $M$, and let $\Gamma$
bound an embedded disk. Does there exist an embedded minimal disk which minimizes the area among the embedded disks
with boundary $\Gamma$?

\vspace{0.3cm}

This is the most general case for a curve to bound an embedded minimal disk. This is because if $\Gamma$ does not bound
any embedded disk in $M$, then, of course, there is no embedded minimal disk bounding $\Gamma$ at all. Our main result
is as follows:

\vspace{0.3cm}

\noindent {\bf Theorem 3.1:} Let $\Gamma$ be a simple closed curve bounding an embedded disk in a closed $3$-manifold
$M$. Then, there exists a disk $\Sigma$ in $M$ with $\partial \Sigma =\Gamma$ such that $\Sigma$ minimizes the area
among all the embedded disks bounding $\Gamma$. Moreover, $\Sigma$ is minimal and smoothly embedded everywhere except
where the boundary $\Gamma$ meets the interior of $\Sigma$.

\vspace{0.3cm}

In particular, if $\Sigma$ is as in the theorem, then there is a continuous parametrization $\varphi: D^2 \rightarrow
M$ of $\Sigma$, with $\varphi(D^2) = \Sigma$ and $\varphi (\partial D^2) = \Gamma$, such that $\varphi$ is smooth
embedding on $D^2- \varphi^{-1}(\Gamma)$ and the image $\Sigma - \Gamma$ is a minimal surface (See Figure 1).
$\lambda=\Gamma\cap\varphi(int(D^2))$ is known as {\em the coincidence set} in the literature (See Remark 3.2). If
$\lambda=\emptyset$, then $\Sigma$ is a smooth embedded minimal disk in $M$ with boundary $\Gamma$. We call such a disk
$\Sigma$ as a {\em pseudo-area minimizing disk}.

The outline of the technique is summarized as follows: Let $\Gamma$ be any simple closed curve bounding an embedded
disk in a manifold $M$. By drilling out a small neighborhood $N_i$ of $\Gamma$ and changing the metric in a very small
neighborhood of the boundary, we can get a manifold $M_i$ with a convex boundary $\partial{M_i}$. For a curve $\Gamma_i
\subset \partial{M_i}$ homotopic to $\Gamma$ in $M$, there is an area minimizing embedded disk $D_i$ in $M_i$ with
$\partial D_i=\Gamma_i$ by \cite{MY2}. When $N_i$ gets smaller, we get a sequence of embedded disks $\{D_i\}$ in $M$
where the areas of the disks are approaching to the minimum area for embedded disks, and $\partial{D_i}=\Gamma_i
\rightarrow \Gamma$. Then, the idea is to obtain a limiting surface out of this sequence and to analyze its structure.

By using the standard generalizations on $M$ to ensure the embeddedness of the solutions of the Plateau problem in
Meeks-Yau setting, we also give a generalization of the main result to a homogeneously regular $3$-manifold $M$ with
sufficiently convex boundary (See Theorem 3.2).

On the other hand, by slightly modifying the pseudo-area minimizing disk $\Sigma$ with boundary $\Gamma$, it is easy to
get a smooth, embedded disk $\Sigma'$ with boundary $\Gamma$ such that $\Sigma'$ is minimal everywhere except for a
very small region (See Corollary 3.3).

The organization of the paper is as follows: In Section 2, we cover some basic results which will be used in the
remaining part of the paper. In section 3, we prove the main result. Then in section 4, we make some final remarks.

\subsection{Acknowledgements:}

I would like to thank Danny Calegari, David Gabai and Joel Hass for very helpful comments.

\begin{figure}[t]

\relabelbox  {\epsfxsize=6in

\hspace{1.5cm} \centerline{\epsfbox{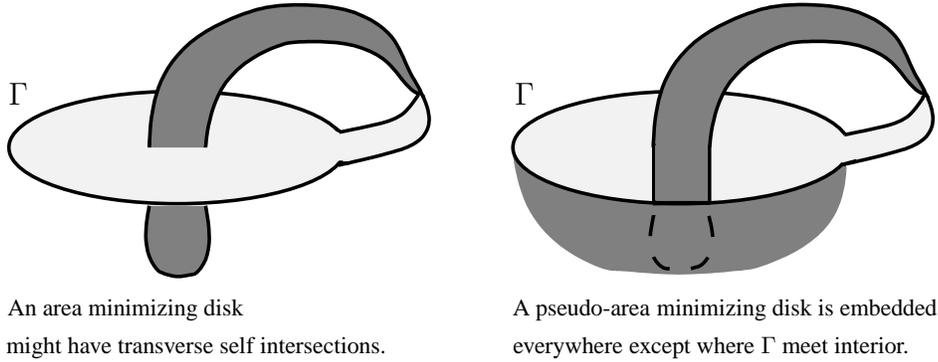}}}

\relabel{1}{$\Gamma$}

\relabel{2}{\scriptsize An area minimizing disk}

\relabel{3}{\scriptsize might have transverse self intersections.}

\relabel{4}{$\Gamma$}

\relabel{5}{\scriptsize A pseudo-area minimizing disk is embedded}

\relabel{6}{\scriptsize everywhere except where $\Gamma$ meet interior.}

\endrelabelbox

\caption{\label{fig:figure1} \small {For a given $\Gamma$, while the area minimizing disks might have transverse self
intersections (left), the pseudo-area minimizing disks have no transverse self-intersection and are embedded everywhere
except where the boundary bumps into interior (right).}}

\end{figure}

\section{Preliminaries}

In this section, we will overview the basic results which we use in the following sections. For more details on these,
see\cite{CM} or \cite{HS}.

\begin{defn} Let $M$ be a $3$-manifold. A {\em minimal disk} in $M$ is a disk whose mean curvature vanishes everyhere.
An {\em area minimizing disk} in $M$ is a disk which has the smallest area among the disks in $M$ with the same
boundary. A {\em pseudo-area minimizing disk} in $M$ is a disk which has the smallest area among the embedded disks in
$M$ with the same boundary, and has no transverse self-intersection.
\end{defn}

\begin{defn} \cite{HS} Let $M$ be a compact Riemannian $3$-manifold with boundary. Then $M$ is {\em mean convex}
(or $M$ has sufficiently convex boundary) if the following conditions hold.

\begin{itemize}

\item $\partial M$ is piecewise smooth.

\item Each smooth subsurface of $\partial M$ has nonnegative curvature with respect to inward normal.

\item There exists a Riemannian manifold $N$ such that $M$ is isometric to a submanifold of $N$ and
each smooth subsurface $S$ of $\partial M$  extends to a smooth embedded surface $S'$ in $N$ such
that $S' \cap M = S$.

\end{itemize}

\end{defn}

\begin{defn} \cite{MY1}, \cite{HS} Let $M$ be a Riemannian $3$-manifold. Then $M$ is {\em homogeneously regular}
if the sectional curvature is bounded above, and injectivity radius is bounded away from $0$.
\end{defn}

In this paper, we will use the following definition for extreme curves. Note that this definition
is different from the one in the literature (lying in the boundary of its convex hull), and our
definition is more general than the other one.

\begin{defn}
$\Gamma\subset M$ is an {\em extreme curve} if it is a curve in the boundary of a mean convex
submanifold in $M$.
\end{defn}

Now, we state the main facts which we use in the following sections.

\begin{lem}\cite{MY1}, \cite{MY2}
Let $M$ be a compact manifold with sufficiently convex boundary, and $\Gamma$ be a nullhomotopic simple closed curve in
$\partial M$. Then, there exists an area minimizing disk $D$ in $M$ with $\partial D = \Gamma$. Moreover, unless
$D\subset \partial M$, all such disks are properly embedded (i.e. the boundary of the disk is in the boundary of the
manifold) in $M$ and they are pairwise disjoint.
\end{lem}

Now, we state two lemmas due to Hass and Scott \cite{HS}, which we use in the following sections.

\begin{lem} (\cite{HS}, Lemma 3.1) Let $M$ be a closed Riemannian $3$-manifold. Then, there exists
an $\epsilon>0$ such that for any $x\in M$, the ball $B_\epsilon(x)$ of radius $\epsilon$ about $x$ in $M$ has the
following property: If  $\Gamma\subset \partial B_\epsilon(x)$ is a simple closed curve, and if $D$ is an area
minimizing disk in $M$ with $\partial D =\Gamma$, then $D$ is properly embedded in $B_\epsilon(x)$.
\end{lem}

\begin{lem} (\cite{HS}, Lemma 3.6) Let $M$ be a compact Riemannian $3$-manifold with strictly convex
boundary. Let $\{D_i\}$ be a sequence of properly embedded area minimizing disks in $M$ which have
uniformly bounded area. Then there is a subsequence $\{D_{i_j}\}$ which converges to a collection
of properly embedded area minimizing disks. If $\{D_i\}$ has a limit point, then the collection is
not empty.
\end{lem}

\section{Main Result}

In this section, we will prove the main result of the paper.

Let $\Gamma$ be a simple closed curve which is the boundary of an embedded disk $E$ in a closed
Riemannian $3$-manifold $M$.

First, we will construct a sequence of embedded {\em almost} area minimizing disks $\{D_i\}$ in $M$ with $\partial D_i
= \Gamma_i$ and $\Gamma_i \rightarrow \Gamma$ by using the techniques of Calegari and Gabai in \cite{Ga} and \cite{CG}.
Then, by taking the limit of this sequence as in \cite{Ga} and \cite{HS}, we will get a limit object $\Delta$. Then, by
analyzing this object, we will show the main result of the paper.

\subsection{The Sequence} \

\vspace{.3cm}

Take a sequence of open solid tori $\{N_i\}$ which are neighborhoods of $\Gamma$. That is, fix a sufficiently small
$\epsilon
>0$, and let $N_i= N_{\frac{\epsilon}{i}}(\Gamma)$, where $N_\epsilon(.)$ represents the open $\epsilon$ neighborhood
in $M$. Then $N_1 \supset N_2 \supset ....\supset N_i \supset N_{i+1} \supset ... $ and $\bigcap_{i=1}^\infty N_i =
\Gamma$.

Now, let $M_i = M - N_i$. Clearly, $\{M_i\}$ are compact $3$-manifolds with torus boundary.
Moreover, $M_1 \subset M_2 \subset ... \subset M_i \subset M_{i+1} \subset ...$ and
$\bigcup_{i=1}^\infty M_i = M-\Gamma$. Also, note that for $\delta<\epsilon$, $\partial
\overline{N_\delta(\Gamma)}$ is a torus, and let $F:(0,\epsilon] \rightarrow \BR$ be a function
such that $F(\delta)= |\partial \overline{N_\delta(\Gamma)}|$ where $|.|$ represents the area in
$M$. Since $\partial \overline{N_\delta(\Gamma)}$ degenerates into $\Gamma$ as $\delta\rightarrow
0$, then $F(\delta)\to 0$ as $\delta \to 0$. Hence, $|\partial M_i| \to 0$ as $i\to \infty$.

Now, we will construct a sequence of area minimizing disks in $M$. Let $E$ be the disk in $M$
bounding $\Gamma$. Modify $E$ if necessary so that $E$ is transverse to $\partial M_i$. Then, let
$\widehat{\Gamma}_i = E\cap \partial M_i$. In other words, for each $i$, let $\widehat{\Gamma}_i$
be a simple closed curve in $\partial M_i= \partial N_i$ which is isotopic to the core curve
$\Gamma$ in the solid torus $N_i$.

Consider the manifolds with torus boundary $M_i \subset M_{i+1}$. Change the metric of $M_{i+1}$ in
$M_{i+1} - M_i$ so that $M_{i+1}$ becomes a compact manifold with sufficiently convex boundary, say
$\widehat{M}_{i+1}$ ($M_{i+1}$ with a new metric). Note that the new metric is the same with
original metric of $M$ in $M_i$ part. As $\widehat{M}_{i+1}$ is mean convex, by Lemma 2.1, there
exists an area minimizing disk $\widehat{D}_{i+1}$ in $\widehat{M}_{i+1}$ with $\partial
\widehat{D}_{i+1} = \widehat{\Gamma}_{i+1}$ (See Figure 2).

\begin{figure}[b]

\relabelbox  {\epsfxsize=5in

\centerline{\epsfbox{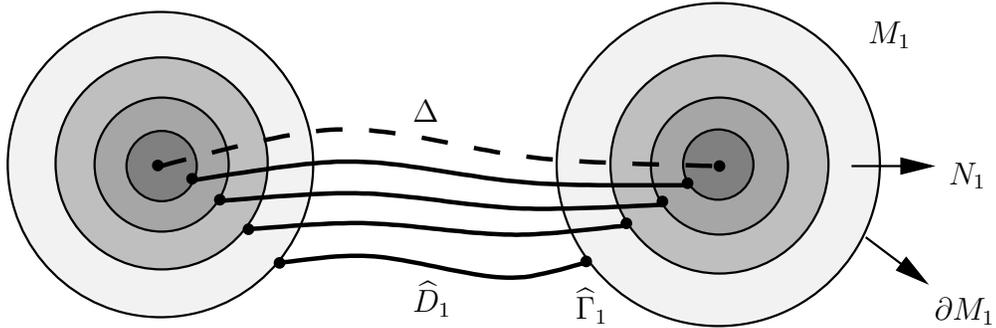}}}

\relabel{1}{$M_1$}

\relabel{2}{$N_1$}

\relabel{3}{$\widehat{D}_1$}

\relabel{4}{$\Delta$}

\relabel{5}{$\widehat{\Gamma}_1$}

\relabel{6}{$\partial M_1$}

\endrelabelbox

\caption{\label{fig:figure2} \small {For any $i$, $M_i = M - N_i$, then $N_i \supset N_{i+1}$, and
$M_i\subset M_{i+1}$. $\widehat{\Gamma}_i \subset \partial M_i$ is a simple closed curve isotopic
to $\Gamma$ in $M$ ($\{\widehat{\Gamma}_i\}$ are shown as pairs of points, and $\Gamma$ is shown as
the innermost pair of points  in the figure). $\widehat{D}_i$ is an area minimizing disk in
$\widehat{M}_i$ ($M_i$ with modified metric) with $\partial \widehat{D}_i = \widehat{\Gamma}_i$.}}

\end{figure}

Now, consider the intersection $\partial M_i \cap \widehat{D}_{i+1}$. By modifying $M_i$ if
necessary, we can assume the intersection is transverse, and it is a collection of simple closed
curves. By construction, the curves in the intersection are either essential in $\partial M_i$ and
isotopic to $\Gamma$ in M, or not essential in $\Gamma$. Let $\{\alpha_1,\alpha_2, ..., \alpha_n\}$
be essential curves, and let $\{\beta_1,\beta_2, ..., \beta_m\}$ be the nonessential ones. By
construction, we know that $n\geq 1$, and $m\geq 0$. Without loss of generality, let $\alpha_1$ be
the innermost curve in $\widehat{D}_{i+1}$ among $\{\alpha_1,\alpha_2, ..., \alpha_n\}$. Let $E_i$
be the subdisk in $\widehat{D}_{i+1}$ with $\partial E_i = \alpha_1$.

Now, consider that $\partial E_i \subset \partial M_i$ and $E_i \subset M_{i+1}$. If $E_i \subset
M_i$, then define the $i^{th}$ element of the desired sequence $D_i = E_i$. Otherwise, $E_i\cap
(M_{i+1} - M_i) \neq \emptyset$, and consists of planar surfaces whose boundary is some curves in
$\{\beta_1,\beta_2, ..., \beta_m\}$. Now replace these planar surfaces in $E_i$ with the isotopic
surface in $\partial M_i$ with the same boundary. Then, smooth out the corners into $M_i$ and push
the parts in $\partial M_i$ into $M_i$ so that surgered $E_i$, say $\widehat{E}_i$, becomes a disk
in $M_i$. Hence, define $\Gamma_i$ as $\alpha_1$, and $D_i$ as $\widehat{E}_i$.

By construction $\{D_i\}$ has the following properties:

\begin{itemize}

\item $D_i$ is a properly embedded smooth disk in $M_i$, i.e. $\partial D_i = \Gamma_i = D_i \cap \partial M_i$.

\item $|D_i| < T_i + 2|\partial M_i|$ where $T_i$ is the infimum of the areas of the embedded disks in $M$ with boundary $\Gamma_i$.

\end{itemize}

By construction, the first property is clear. To see the second property, first consider that $E_i$
is the area minimizing disk in $\widehat{M}_{i+1}$, and the metric of $\widehat{M}_{i+1}$ coincides
with the original metric on $M$ in $M_i$. Clearly, $T_i' < T_i + |\partial M_i|$ where $T_i'$ is
the infimum of the areas of the embedded disks in $M_i$ with boundary $\Gamma_i$. As $D_i$ obtained
by modifying $E_i\cap M_i$, by patching the missing parts from $\partial M_i$ and smoothing out
corners. By \cite{MY2}, smoothing out the folding curves decreases the area, hence $|D_i| \leq
|E_i\cap M_i| + |\partial M_i|$. Since $|E_i\cap M_i| \leq T_i'$, the second property follows.

It is clear that $D_i$ may not be area minimizing in $M$, however, with the second property, they
can be thought as {\em almost area minimizing} in $M$. Note that $E_i\cap M_i\subset D_i$ is area
minimizing in $M$, so only the patched parts of $D_i$ (replacements of $E_i\cap (M_{i+1} - M_i)$)
are not area minimizing, which are a very small part of $D_i$. This is because it is reasonable to
think $D_i$ as almost area minimizing.

Now, let $T$ be the infimum of the areas of the embedded disks in $M$ with boundary $\Gamma$. Let
$A_i$ be the infimum of the areas of embedded annuli with boundary $\Gamma \cup \Gamma_i$. By
construction, $T_i \leq T + A_i$ and $T\leq T_i + A_i$. Since $\Gamma_i \to \Gamma$, $A_i\to 0$ as
$i\to \infty$. This implies $T_i\to T$ as $i\to \infty$. Hence, by the second property and
$|\partial M_i|\to 0$, $|D_i| \to T$ as $i\to \infty$.

To sum up, we constructed a sequence of embedded smooth disks $\{D_i\}$ in $M$ such that $\partial D_i = \Gamma_i
\subset \partial M_i$ with $\Gamma_i \to \Gamma$ and $|D_i| \to T$ as $i\to \infty$
where $T$ is the infimum of the areas of the embedded disks in $M$ with boundary $\Gamma$.\\

\subsection{The Limit} \ \\

In this section, by using the techniques of \cite{Ga} and \cite{HS}, from the sequence $\{D_{i}\}$ of {\em almost} area
minimizing embedded disks in a $M$, we will get a limit object $\Delta$ which can be thought as a special topological
limit. Then, by using this object, we will prove the main result of the paper. The following is a modified version of a
definition of special topological limit due to Gabai \cite{Ga}.

\begin{defn} A collection of pairwise disjoint embedded surfaces $\Delta$ in a Riemannian manifold $M$ is called the {\em topological limit} of the sequence $\{D_{i}\}$
 of embedded surfaces in $M$ if

\begin{itemize}

\item $\Delta = \{ \ x=\lim x_i \ | \ \ x_i \in D_i , \{x_i\} \mbox{ is a convergent sequence in } M \
\}$

\item $\Delta = \{ \ x=\lim x_{n_i} \ | \ \ x_i \in D_i , \{x_i\} \mbox{ has a convergent subsequence }
\{x_{n_i}\} \mbox{ in } M \ \}$

\end{itemize}
\end{defn}

In other words, the sequence $\{D_i\}$ is such that the set of the limits of all $\{x_i\}$ with
$x_i\in M$ and the set of the limits of the subsequences are the same. This is an essential
condition on $\Delta$ to be a collection of \textit{pairwise disjoint} embedded surfaces.
Otherwise, one might simply take a sequence such that $D_{2i+1} = \Sigma_1$ and $D_{2i} = \Sigma_2$
where $\Sigma_1$ and $\Sigma_2$ are intersecting disks. Then, without the first condition ($\Delta$
being just the union of limit points), $\Delta= \Sigma_1 \cup \Sigma_2$ in this case, which is not
a collection of pairwise disjoint embedded disks. However, the first condition forces $\Delta$ to
be either $\Sigma_1$ or $\Sigma_2$, not the union of them. By similar reasons, this condition is
also important to make sure the embeddedness of the disks in the collection $\Delta$.

Now, we will show that there is a subsequence of the sequence constructed in previous part which gives a topological
limit $\Delta$. Then, by showing that the limit $\Delta$ is a collection of disks, we will prove the main result of the
paper.

\begin{thm}
Let $\Gamma$ be a simple closed curve bounding an embedded disk in a closed $3$-manifold $M$. Then, there exists a disk
$\Sigma$ in $M$ with $\partial \Sigma =\Gamma$ such that $\Sigma$ minimizes the area among the embedded disks bounding
$\Gamma$. Moreover, $\Sigma$ is minimal and smoothly embedded everywhere except where the boundary $\Gamma$ meets the
interior of $\Sigma$.
\end{thm}

\begin{pf} Let $\{D_i\}$ be as defined in the previous part. We abuse the notation by taking $\{D_i\}$
as a sequence of open disks. Our aim is to get a convergence as in Definition 3.1 for an
appropriate subsequence of $\{D_i\}$.\\

\noindent {\bf Step 0:} After passing to a subsequence of $\{D_i\}$, the following holds:\\

$\Delta = \{ \ x=\lim x_i \ | \ \ x_i \in D_i , \{x_i\} \mbox{ is a convergent sequence in } M \
\}$

\hspace{.4cm} $= \{ \ x=\lim x_{n_i} \ | \ \ x_i \in D_i , \{x_i\} \mbox{ has a convergent
subsequence } \{x_{n_i}\} \mbox{ in } M \ \}$\\

\begin{pf}
For each $j$ subdivide $M$ into a finite number of closed regions such that the $j+1^{st}$
subdivision is a subdivision of the $j^{th}$ one. Also, let the mesh of these subdivisions
converges to $0$. In other words, let $B^j_k$ be the $k^{th}$ region of $j^{th}$ subdivision and
$M=\bigcup_{k=1}^{n_j}B^j_k$ where $B^{j-1}_i=B^j_{i_1} \bigcup$...$\bigcup B^j_{i_r}$. Also, for
any $j,k$, $diam(B^j_k) < C_j$ where $C_j \rightarrow 0$ as $j\rightarrow \infty$.

Now, choose a subsequence of $\{D_i\}$ such that if $i\geq j$ and $D_i\bigcap B^j_k\neq\emptyset$,
then for any $k>i$, $D_k\bigcap B^j_r \neq \emptyset$. By abuse of notation, replace this
subsequence with the original sequence $\{D_i\}$. Now, for this new sequence, for any $x=\lim
x_{n_i}$ where $x_{n_i}\in D_{n_i}$, by construction, there is a convergent sequence $\{x_i\}$ with
$x_i\in D_i$ such that $x=\lim x_i$. Hence, Step 0 follows.

\end{pf}

\noindent {\bf Step 1:} $\Delta$ is not empty.\\

\begin{pf}
In order to show that $\Delta$ is nonempty, it will suffice to construct a convergent sequence
$\{x_{n_i}\}$ in $M$ with $x_{n_i} \in D_{n_i}$. Now, consider a meridian curve $\gamma'$ of
$\partial M_1$ (Think of $\Gamma_1$ as the longitude of $\partial M_1$). Push $\gamma'$ into $M_1$
a little bit (a small isotopy), and call the new curve as $\gamma$ in $M_1$. Now, $\gamma$ links
$\Gamma$ and $\Gamma_i$ in $M$ for all $i$. Hence, for any $i$, $D_i \cap \gamma$ is not empty. Let
$x_i \in D_i \cap \gamma$. Since $\gamma$ is a simple closed curve, the sequence $\{x_i\}$ must
have a convergent subsequence. Hence, $\Delta$ is nonempty.

\end{pf}

\noindent {\bf Step 2:} Let $Z=\{z_i\}$ be a countable dense subset of $\Delta$ where $Z\cap\Gamma
=\emptyset$. Then, after passing to a subsequence of $\{D_j\}$ the following holds. For any $i$,
there exists a sequence $\{E^i_j\}$ of embedded disks $E^i_j\subset D_j$ which converges to a
smoothly embedded disk $E_i$ such that $z_i\in E_i$.\\

\begin{pf}
As $M$ is a closed manifold, by Lemma 2.2, there exists a $\rho>0$ such that for any $x \in M$, if
$\Gamma \subset \partial B_\rho(x)$, then any area minimizing disk $D$ in $M$ with boundary
$\Gamma$ must be in $B_\rho(x)$.

Now, let $Z_1 = \{z_i\in Z \ | \ d(z_i,\Gamma) < 2\rho \}$ and $Z_2= Z-Z_1$. Also, let $\rho_i =
\frac{d(z_i,\Gamma)}{2}$ for any $z_i \in Z_1$. Now, we claim that for any $z_i \in Z_1$, there is
an embedded disk $E_i$ in $\Delta$ with $z_i\subset E_i$ and $\partial \overline{E}_i \cap
B_{\rho_i}(z_i)= \emptyset$. We also claim that for $z_i\in Z_2$, there is an embedded disk $E_i$
in $\Delta$ with $z_i\subset E_i$ and $\partial \overline{E}_i \cap B_\rho(z_i)= \emptyset$.

Now, fix $z_i \in Z_2$. Since $z_i\in \Delta$, there exists a sequence $\{x^i_j\}$ with $x^i_j\in
D_j$ and $x^i_j\to z_i$. By deleting the tangential points in the intersection, and by modifying
$\rho$ if necessary, we can assume that $\partial B_\rho(z_i)$ is transverse to $\{D_j\}$. Let
$E^i_j= B_\rho (z_i) \cap D_j$. For sufficiently large $j$, being away from $\partial M_j$, almost
area minimizing disk $D_j$ is area minimizing near $z_i$ and hence $E^i_j$ is area minimizing disk
in $B_\rho(z_i)$ by Lemma 2.2. Hence $ |E^i_j| < \frac{1}{2} | \partial B_\rho(z_i)|$ for any $j$,
and the sequence of properly embedded area minimizing disks in $B_\rho(z_i)$ have uniformly bounded
area. Also, since $M$ is closed, we can assume $B_\rho(z_i)$ has strictly convex boundary, by
taking a smaller $\rho$ if necessary (Since $M$ is compact, there exists $R$ such that for any
$0<r<R$, $B_r(x)$ has strictly convex boundary for any $x\in M$). Hence, by Lemma 2.3, $\{E^i_j\}$
has a subsequence converging to $E_i$ where $E_i$ is the area minimizing disk in $M$ with $z_i\in
E_i$. By using diagonal subsequence argument, we can modify our sequence $\{D_i\}$ accordingly. In
other words, for each $z_i$, use above argument and get a diagonal subsequence, and call this new
sequence $\{D_i\}$ again (abuse of notation), and define $\Delta$ for this new sequence.

Similarly, fix $z_i\in Z_1$. As above, let $\{x^i_j\}$ be a sequence with $x^i_j\in D_j$ and
$x^i_j\to z_i$, and $E^i_j= B_{\rho_i} (z_i) \cap D_j$. By construction, $\rho_i\leq \rho$ for any
$i$. Again for sufficiently large $j$, being away from $\partial M_j$, almost area minimizing disk
$D_j$ is area minimizing near $z_i$ and hence $E^i_j$ is area minimizing disk in $B_{\rho_i}(z_i)$
by Lemma 2.2. As before, by Lemma 2.3, $\{E^i_j\}$ has a subsequence converging to $E_i$ where
$E_i$ is the area minimizing disk in $M$ with $z_i\in E_i$. Continuing with the diagonal
subsequence argument for each $z_i$, we get a subsequence of $\{D_i\}$ with the required
properties.
\end{pf}

Hence, with above construction, we get $\bigcup_{z_i\in Z} E_i = \Delta - \Gamma$. Now, we will analyze the structure
of the limit object $\Delta$. In particular, we will prove that $\Delta$ is a union of ``disks" which minimize area
among the embedded disks. Now, clearly $\Gamma \subset \Delta$. We want to specify some parts of $\Gamma$ which bumps
into interior of $\Delta$. In other words, define {\em the coincidence set} $\lambda\subset\Gamma$ such that

$$\lambda= \{ \ x\in\Gamma \ | \ \exists \rho_x>0 , \exists \{y_j\} \mbox{ with } y_j\in D_j \mbox{ and }
\hat{d}_j(y_j,\Gamma_j)>\rho_x \mbox{ such that } y_j\to x\}$$

Here, $\hat{d}_j$ is the induced path metric on $D_j$. In other words, a point in the boundary $\Gamma$ is in the
coincidence set $\lambda$, if there exists a sequence $\{z_j\}$ in $\Delta$ which is away from the boundary $\Gamma$
and $z_j\to x$. Hence by definition, the coincidence set $\lambda$ corresponds to some closed subsegments of $\Gamma$
which meets the interior part of $\Delta$. It can be thought as the {\em defective} parts of the embedded minimal disks
we are constructing.
Note that $\lambda$ might be empty.\\

\noindent {\bf Step 3:} $\Delta-\Gamma$ is a minimal surface in $M$.\\

\begin{pf}
By Step 2, for each $z_i\in Z$, there is an embedded disk $E_i \subset \Delta$. Now, if $x\in
E_i\cap E_j$, then they must coincide in a neighborhood of $x$. Otherwise, since $E_i$ and $E_j$
are minimal disks, they must cross transversely near $x$ \cite{HS}. However, by construction, this
would imply that $D_i$ is not embedded for sufficiently large $i$. Hence, $E_i$'s are either
pairwise disjoint or locally coincide. As $\{z_i\}$ is a dense subset of $\Delta$, by Step 2, for
any $x \in \Delta - \Gamma$, we can find a neighborhood of $x$ in $\Delta$, say $E_x$, which is an
open minimal disk. Minimality comes from being locally area minimizing. Hence, $\Delta-\Gamma =
\bigcup_{x\in \Delta-\Gamma} E_x$. This shows that each component $\Sigma_i '$ of $\Delta - \Gamma$
is a surface.
\end{pf}

\noindent {\bf Step 4:} $\Delta$ is a collection of disks which minimize area among the embedded disks.\\

\begin{pf}
Now, we will show that for each component $\Sigma_i '$ of $\Delta - \Gamma$, $\Sigma_i = \Sigma_i'\cup \Gamma$ is a
pseudo-area minimizing disk in $M$ where $\Sigma_i - \Gamma$ is a smooth minimal surface. In other words, we will show
that there is a continuous map $\varphi_i:D^2 \to M$ such that $\varphi(\partial D^2) = \Gamma$ and $\varphi (D^2) =
\Sigma_i$ with $\varphi$ is an embedding except at $\varphi^{-1}(\lambda)$. i.e. $\Sigma_i\subset M$ is the image of a
continuous map from a disk where it fails to be an embedding only at the coincidence set $\lambda$.

Let $\alpha$ be a simple closed curve in $\Sigma_i '$. Let $A$ be a neighborhood of $\alpha$ which is a very thin
annulus in $\Sigma_i '$. Now, let $g:D\to A$ be an isometric immersion of a disk $D$ into $\Sigma_i$ where $D$ is a
very long thin rectangle with $g(D)=A$. Also, assume that $|D|> C$ where $C$ is a constant with $C> |D_i|$ for any $i$.
Existence of such a $C$ comes from the construction as $D_i$ is area minimizing disk in $M_i$, and $|D_i| < |A_i| +
|D_1|$ where $A_i$ is a very thin annulus with $\partial A_i = \Gamma_1\cup\Gamma_i$. Since we can find a uniform bound
for $|A_i|$, the existence of $C$ follows.

Now, as $D_i \to \Delta$, and $A \subset \Delta - \Gamma$, we can find isometric immersions $g_i: D^2 \to D_i$ such
that $g_i \to g$ in $C^\infty$ topology. Now, there are two cases. Either $g_i(D^2)\subset D_i$ is also a thin annulus
$A_i$ in $D_i$ approximating $A$, or $g_i(D^2)$ is an embedded disk in $D_i$ which is spiraling around $A$. In the
latter case, it would mean that $|D_i|>|D|> C$ which is a contradiction. Hence, $g_i(D^2) = A_i$. Then, we can choose a
suitable essential simple closed curve $\beta_i$ in each annulus $A_i$ ($\beta_i$ is a core curve of the annulus $A_i$)
such that the sequence $\{\beta_i\}$ converges to $\alpha$, i.e. $\beta_i \to \alpha$. Let $F_i$ be the disks in $D_i$
with $\partial F_i = \beta_i$.

Now, further assume that $\alpha$ separates $\Sigma_i$ into two parts, say $S_1$ and $S_2$, and $\Gamma \subset S_2$.
In other words, $S_1\subset \Sigma_i$ and $\partial S_1 = \alpha$ with $S_1\cap \lambda = \emptyset$. Being away from
$\Gamma$, hence from the coincidence set $\lambda$, this implies that the disks $\{F_i\}$ are area minimizing disks
(not almost area minimizing). Hence, as in Lemma 3.3 of \cite{Ga}, the sequence of disks $\{F_i\}$ converges to a disk
$\Omega$ in $\Sigma_i$ with $\partial \Omega = \alpha$, i.e. $S_1=\Omega$ is a smooth minimal disk in $M$. Hence, we
show that for any separating simple closed curve $\alpha$ in $\Sigma_i$ with $S_1\cap \lambda = \emptyset$, there is a
smooth disk $\Omega$ in $\Sigma_i '$ with $\partial \Omega =\alpha$.

Now, by choosing a suitable sequence of simple closed curves $\{\alpha_n\}$ in $\Sigma_i$, we can exhaust $\Sigma_i$
with disks $\Omega_n$ such that $\Omega_1 \subset \Omega_2 \subset ... \subset \Omega_n \subset ...$ with $\Sigma_i =
\bigcup_n \Omega_n$ where $\Omega_n$ is a disk in $\Sigma_i$ with $\partial \Omega_n = \alpha_n$ and $\Omega_n \cap
\lambda =\emptyset$. This means that $\alpha_n \to \Gamma\cup l\cup \lambda$ where $l$ is a collection of line segments
$\{l^k_j\}$ in $\Sigma_i$ which connects $\Gamma$ with one of the endpoints of the components $\{\lambda_k\}$ of
$\lambda$. In particular, for each component $\lambda_k$ of $\lambda$, let $n_k$ be the number of line segments
$\{l^k_1,l^k_2,..,l^k_{n_k}\}$ connecting the component $\lambda_k$ to the $\Gamma$. Then $n_k$ is the number of the
{\em local sheets} of $\Sigma_i$ near $\lambda_k$. i.e. the components of $N(\lambda_k) \cap \Omega_n$ where
$N(\lambda_k)$ is a sufficiently small neighborhood of $\lambda_k$ in $M$ and $n$ is sufficiently large.

Since the sequence of the disks $\{\Omega_n\}$ exhausting $\Sigma_i$ does not contain $\lambda$, $\lambda$ is in the
component $\Sigma_i - \Omega_n$ which also contains $\Gamma$. Hence, for sufficiently large $n$, $\Sigma_i - \Omega_n$
is the union of some thin neighborhood of $\Gamma$ and some thin strips around line segments $\{l^k_j\}$ connecting
$\Gamma$ to $\lambda$, as $\Omega_n$ is exhausting $\Sigma_i$. Hence, $\alpha_n$ curves approach $l\cup\lambda$ from
both sides and $\partial \Omega_n = \alpha_n \to \Gamma\cup l\cup\lambda$. Hence, we can get a continuous
parametrization $\hat{\varphi}: D^2\to M$ with $\hat{\varphi}(D^2)=\Sigma_i$ and $\hat{\varphi}(\partial D^2) =
\Gamma\cup l\cup\lambda$.

Now, our aim is to get a continuous map $\varphi:D^2\to M$ with $\varphi(D^2) = \Sigma_i$ and $\varphi(\partial
D^2)=\Gamma$ by modifying $\hat{\varphi}$. Fix $k_0$ and $j_0$. Consider $l^{k_0}_{j_0}\cup \lambda_{k_0}$. Let
$U=N(l^{k_0}_{j_0}\cup \lambda_{k_0})$ be a very small neighborhood of $l^{k_0}_{j_0}\cup \lambda_{k_0}$ in $\Sigma_i$
such that $U \cap l^k_j = \emptyset$ for any $(k,j)\neq (k_0,j_0)$. Consider $\hat{\varphi}^{-1}(U) \subset D^2$. Let
$V$ be the component of $\hat{\varphi}^{-1}(U)$ in $D^2$ which contains the segment
$\hat{\varphi}^{-1}(l^{k_0}_{j_0})\subset \partial D^2$. Hence, $\hat{\varphi}|_V:V\to U$ and $\hat{\varphi} (V\cap
\partial D^2) \subset \Gamma \cup \l^{k_0}_{j_0} \cup \lambda_k$.

Now, we can continuously modify $\hat{\varphi}$ in $V$ into a continuous map $\varphi$ so that $\varphi|_V:V\to U$ is a
continuous embedding with $\varphi (\partial\overline{V}\cap int(D^2)) = \hat{\varphi} (\partial\overline{V}\cap
int(D^2))$ and $\varphi(V\cap \partial D^2)\subset \Gamma$. To see this, one can define a continuous family of maps
$\hat{\varphi}_t:V\to U$ with $0\leq t \leq 1$ such that $\hat{\varphi}_0=\hat{\varphi}|_V$ and
$\hat{\varphi}_t(V\cap\partial D^2)\subset \Gamma\cup (l^{k_0}_{j_0}\cup \lambda_{k_0})_t$. Here, $(l^{k_0}_{j_0}\cup
\lambda_{k_0})_t$ is a subsegment in $l^{k_0}_{j_0}\cup \lambda_{k_0}$ which is getting smaller as $t$ increases, and
finally $(l^{k_0}_{j_0}\cup \lambda_{k_0})_1$ is the endpoint of \ $l^{k_0}_{j_0}\cup \lambda_{k_0}$ in $\Gamma$.
Intuitively, this continuous deformation of the parametrizations corresponds to pushing $l^{k_0}_{j_0}\cup \lambda_k$
into $\Gamma$ in $\Sigma_i$. Since, $U=N(l^{k_0}_{j_0}\cup \lambda_{k_0})$ is disjoint from other line segments
$l^k_j$, one can modify $\hat{\varphi}$ for each $l^k_j$. Finally, we get a continuous map $\varphi:D^2\to M$ with
$\varphi(D^2) = \Sigma_i$ and $\varphi(\partial D^2)=\Gamma$.

Even though we obtained $\varphi$ from $\hat{\varphi}$ which may not be smooth along $l\cup \lambda$, since the choice
of $l$ is arbitrary, and $\Sigma_i-\Gamma$ is smooth by previous steps, $\varphi$ can be chosen as a smooth embedding
on $D^2-\varphi^{-1}(\Gamma)$ by construction. If the coincidence set $\lambda=\Gamma \cap \varphi(int(D^2))$ is empty,
then $\Sigma_i$ is an embedded minimal disk in $M$ with $\partial \Sigma_i =\Gamma$. Otherwise, the disk $\Sigma_i$
might fail to be smooth on $\lambda$. Note also that the restriction of $\varphi$ to the interior of $D^2$ may not be
an embedding either, if $n_k>1$ for some $k$ (See Section 4).

Now, we claim that $\Sigma_i$ is minimizing area among the embedded disks with boundary $\Gamma$. Otherwise, there is a
compact subdisk $E$ of $\Sigma_i$ which is not area minimizing among the embedded disks. Then there is an embedded disk
$E'$ in $M$ with smaller area. Let $|E|-|E'|= \xi$. By construction, since $D_i\to\Delta$, we can find sufficiently
close disks in the sequence $\{D_i\}$ to $\Sigma_i$ such that there is a subdisk $E'' \subset D_i$ with $| \ |E| -
|E''| \ | < \frac{\xi}{2}$ and $|\mathcal{A}| < \frac{\xi}{2}$ where $\mathcal{A}$ is an annulus with boundary
$\partial E \cup \partial E''$. However, this implies $|E'| + |\mathcal{A}| < |E''|$. This is a contradiction as $E'
\subset D_i$ is area minimizing among the embedded disks with same boundary. So, $\Sigma_i$ is also area minimizing
among the embedded disks with same boundary. Hence, Step 4 follows.
\end{pf}

Hence, by taking one of the $\Sigma_i$ piece in $\Delta$ given by Step 4, we get a continuous map
$\varphi: D^2 \rightarrow M$ with $\varphi(D^2) = {\Sigma}_i$ and $\varphi (\partial D^2) =
\Gamma$. Moreover, $\varphi$ is a smooth embedding on $D^2-\varphi^{-1}(\Gamma)$, and $\Sigma_i -
\Gamma$ is an embedded minimal surface. The proof follows.
\end{pf}

\begin{rmk}
This theorem shows that for a given simple closed curve $\Gamma$ bounding a disk in a closed $3$-manifold, there exists
a pseudo-area minimizing disk $\Sigma$ in $M$ with $\partial \Sigma = \Gamma$. Moreover, the theorem gives the
structure of the pseudo-area minimizing disks: A pseudo-area minimizing disk $\Sigma$ may not be an embedded disk, but
it can only fail embeddedness if the boundary $\Gamma$ bumps into the interior of $\Sigma$ ($\varphi(int(D^2)$). Also,
$\Gamma$ can only intersect the interior of $\Sigma$ nontransversely, i.e. $\Sigma$ has no transverse
self-intersection. A pseudo-area minimizing disk $\Sigma$ is smooth and minimal everywhere except where $\Gamma$ meets
interior of $\Sigma$.
\end{rmk}

On the other hand, this result is true for more general manifolds, namely homogeneously regular
$3$-manifolds with sufficiently convex boundary. In Section $6$ of \cite{HS}, one might find the
reasons why we need the conditions of being homogeneously regular and being sufficiently convex to
extend these results.

\begin{thm} Let $\Gamma$ be a simple closed curve bounding an embedded disk in a
homogeneously regular $3$-manifold $M$ with sufficiently convex boundary. Then, there exists a disk
$\Sigma$ in $M$ with $\partial \Sigma =\Gamma$ such that $\Sigma$ minimizes the area among the
embedded disks bounding $\Gamma$. Moreover, $\Sigma$ is embedded in the interior, and it is smooth
and minimal everywhere except where the boundary $\Gamma$ meets the interior of $\Sigma$.
\end{thm}

\begin{pf} Since Lemma 2.1 is valid for this new ambient manifolds, we can still construct the
sequence $\{D_i\}$ as before. By replacing the lemmas we used in the previous theorem with their
analogs in the new ambient manifold as in Section $6$ of \cite{HS}, the same proof would work.
\end{pf}

In other words, this theorem applies to closed $3$-manifolds, compact $3$-manifolds with
sufficiently convex boundary, homogeneously regular non-compact $3$-manifolds, and homogeneously
regular non-compact $3$-manifolds with sufficiently convex boundary.

The structure of the pseudo-area minimizing disks given by the main theorem also tells us how to construct nearby
special smoothly embedded disks with the same boundary.

\begin{cor}
Let $M$ be a closed $3$-manifold or a homogeneously regular $3$-manifold with sufficiently convex boundary. Let
$\Gamma$ be a simple closed curve bounding an embedded disk in $M$. Then, for any given $\epsilon>0$, there exist a
smooth, embedded disk $\Sigma_\epsilon$ in $M$ with $\partial \Sigma_\epsilon = \Gamma$ such that $|\Sigma_\epsilon| <
C_\Gamma +\epsilon$ and there exists a small region $R$ in $\Sigma_\epsilon$ where $|R|<\epsilon$ with
$\Sigma_\epsilon-R$ is a minimal surface. Here $|.|$ represents the area, and $C_\Gamma$ is the infimum of the areas of
the embedded disks in $M$ with boundary $\Gamma$.
\end{cor}

\begin{pf}
According to Theorem 3.1, the pseudo-area minimizing disk $\Sigma$ can only fail embeddedness along the coincidence set
$\lambda$ where the boundary $\Gamma$ bumps into the interior of $\Sigma$. Hence, for a given $\epsilon$, take a small
neighborhood of this segment in $\Sigma$, and push every sheet slightly away from $\lambda$ so that we get a smooth,
embedded disk. Since we can choose the neighborhood as small as we want, the corollary follows.
\end{pf}

\begin{rmk} {\bf (Regularity near the coincidence set and the thin obstacle problem)} Another interesting question on the
problem is the regularity of $\varphi$ near the coincidence set. This is a well-known problem in the literature and
it's called as {\em Thin obstacle problem} or {\em Signorini problem}. Clearly, $\varphi$ may not be smooth along the
coincidence set, but it may have some regularity when restricted to one side of the coincidence set. In the classical
setting of the problem in dimension $2$, let $\Omega$ be a bounded open subset in $\BR^2$ and $A$ be a line segment in
$\Omega$. Let $\psi:A\to \BR$ and $g:\partial \Omega \to \BR$ be given where $g \geq \psi$ on $\partial \Omega \cap A$
where $\psi, g$ are smooth. Then the question is to find $\varphi: \Omega \to \BR$ where $\varphi|_{\partial \Omega} =
g$ and $\varphi|_A \geq \psi$, and $\varphi(\Omega)$ has minimum area. In the classical works \cite{Fr}, \cite{Ni} on
the problem, the authors showed the Lipschitz continuity of the solutions in all dimensions, and in dimension $2$,
Richardson showed that the solutions are $C^{1,\frac{1}{2}}$ which is the optimal regularity \cite{Ri}. Recently,
Guillen generalized Richardson's result to any dimension \cite{Gui}.

Hence, when we apply to these results to our case, when $M=\BR^3$, and given simple closed curve $\Gamma \subset \BR^3$
is smooth, then above results imply that if $\varphi$ is as in the main theorem, then for each local sheet of $\Sigma$,
$\varphi$ is $C^{1,\frac{1}{2}}$ on either side of the coincidence set $\lambda$. Also, when $M$ is a hyperbolic
$3$-manifold and $\Gamma$ is a geodesic, then Calegari-Gabai's result when applied to our case imply that $C^1$
regularity of $\varphi$ on either side of the coincidence set (Lemma 1.31 in \cite{CG}).

On the other hand, assuming $\Gamma$ is smooth enough, it is possible to show that $\varphi$ is $H^{1,2}$ near
coincidence set by arguing like in \cite{CG} (Section 1.7). In particular, if $\Gamma$ is a smooth simple closed curve
in $M$, then it can be showed that $\varphi:D^2 \to M$ given in the main theorem is in the Sobolev space $H^{1,2}$,
i.e. the derivative $d\varphi$ is defined and $L^2$ in the sense of distribution. The idea is basically same. If
$\varphi_i: D^2 \to M$ with $\varphi(D^2) = D_i$ where $\{D_i\}$ is the sequence of area minimizing disks constructed
in section 3.1 with modified metric. Then each $\varphi_i$ induce a conformal structure on $D^2$. The $L^2$ norm of the
derivatives can only blow up along a neck pinch. In our case, we work with $D^2$, and hence we cannot have a neck pinch
by energy minimizing property of $\varphi_i$. So, we get the limit of these conformal structures, and $L^2$ norm of $d
\varphi$ can be bounded in terms of the $L^2$ norms of $d \varphi_i$. In other words, the idea is same with the
argument in \cite{CG} except we work with disks instead of surfaces.

There are also relevant results in the literature about the structure of the coincidence set. When $M=\BR^3$ and
$\Gamma$ is analytic, then the coincidence set is a finite union of points and intervals in certain cases \cite{Le},
\cite{At} (See Remark 1.32 in \cite{CG}).
\end{rmk}

\section{Concluding Remarks}

In this paper, we studied the embeddedness question of minimal disks in $3$-manifolds. Unlike the many results
considering the question that ``for which curves must the area minimizing disks be embedded?" in the literature, we
analyzed the structure of the surface which minimizes area among the embedded disks whose boundary is any given simple
closed curve. Hence, we showed that among all the embedded disks with fixed boundary, the area minimizer (pseudo-area
minimizing disk) exists, however it may not be an embedded disk. It is a disk in the manifold, and it only fails to be
an embedding along the coincidence set where the boundary bumps nontransversely into the interior. Other than this
exceptional part, the disk is a smoothly embedded minimal surface. Hence, for any simple closed curve in a $3$-manifold
$M$, we construct a canonical {\em almost} embedded disk in $M$ among the embedded disks bounding the given curve.

Intuitively, one might think the pseudo-area minimizing disk $\Sigma$ in the following way. Let $D$ be the area
minimizing disk $M$ with boundary $\Gamma$. As in Figure 1, $D$ may not be embedded, and it might have self
intersection. Then, one can get $\Sigma$ from $D$ by pushing the self intersection into the boundary. The interesting
fact here is that $\Sigma$ is smooth, minimal, and embedded everywhere except where $\Gamma$ meets interior of
$\Sigma$, say $\lambda$. An alternative way to see the picture is that if you push from the coincidence set $\lambda$
in $\Sigma$ into the convex part, you can reduce the area of $\Sigma$ (like a folding curve in \cite{MY2}), but you
create a transverse self intersection in the interior. This cannot happen as $\Sigma$ is minimizing area among the {\em
embedded} disks. In other words, the coincidence set $\lambda$ behaves like a barrier to embeddedness, even though you
can reduce the area by going in that direction.

Also, one might ask whether $\varphi$ is an embedding on whole $int(D^2)$ or not. This is not true in general. The
reason for that the interior might have nontransverse self-intersection with itself at the coincidence set $\lambda$.
For example, if $\Sigma$ has more than one local sheet near the coincidence set $\lambda$, $\varphi|_{int(D^2)}$ can
not be an embedding. To construct such an example, one can take two ``parallel" close embedded area minimizing disks
$\Sigma_1$ and $\Sigma_2$. Then, by adding a tiny bridge $\beta$ between them one can get another area minimizing disk
$\Sigma$ with new boundary which is close to $\Sigma_1\cup\Sigma_2\cup\beta$. Now, if you make a thin and long horn
from the part close to $\Sigma_1$ which intersect both $\Sigma_1$ and $\Sigma_2$ transversely like in Figure 1, then
the area minimizing disk $\hat{\Sigma}$ which minimizes area for the new boundary curve is an example for such a
situation. One needs to push the intersection of the horn with both $\Sigma_1$ and $\Sigma_2$ to the boundary, and the
interiors of them will meet in the boundary. Hence, $\varphi$ cannot be an embedding on $int(D^2)$ in general.

On the other hand, there is a relevant result about the same question due to Meeks-Yau \cite{MY3}. They give a
necessary condition for a sufficiently smooth simple closed curve in a $3$-manifold to bound a strictly stable embedded
minimal disk.

\begin{thm} (\cite{MY3}, Theorem 3) Let $\Gamma$ be a $C^{2,\alpha}$ simple closed curve in a $3$-manifold $M$.
Then, $\Gamma$ bounds a strictly stable minimal disk $\Sigma$ in $M$ if and only if there exists a
sufficiently convex codimension-$0$ submanifold $N$ in $M$, whose topological type is a $3$-ball,
and $\Gamma \subset \partial N$.
\end{thm}

Hence, by combining this result with ours for $C^{2,\alpha}$ smooth simple closed curves, we conclude that if $\Gamma$
is a $C^{2,\alpha}$ simple closed curve in a $3$-manifold $M$, and there is no sufficiently convex domain $N$ as in
above theorem with $\Gamma \subset \partial N$, then the pseudo-area minimizing disk $\Sigma$ given by our main theorem
is not embedded up to the boundary. In other words, if $\varphi: D^2 \to M$ parameterizes $\Sigma$ with $\varphi (D^2)
= \Sigma$, $\varphi (\partial D^2) = \Gamma$, then $\Gamma \cap \varphi|_{int(D^2)} \neq \emptyset$. However, neither
our result nor the above theorem of Meeks-Yau says anything about unstable minimal disks. It is still possible for such
a simple closed curve in a $3$-manifold $M$ to bound an embedded unstable minimal disk in $M$.

\end{document}